\documentclass[reqno,draft]{amsart}
\usepackage{multicol, color}

\newcommand{\rem}[1]{}

\theoremstyle{plain}
\newtheorem{lemma}{Lemma}
\newtheorem{theorem}[lemma]{Theorem}
\newtheorem{corollary}[lemma]{Corollary}

\newtheorem{proposition}[lemma]{Proposition}

\theoremstyle{remark}

\newcommand*  {\R} {{\mathbb R}}

\newcommand*{\field}[1]{\mathbb{#1}}

\def\aa{\alpha}
\def\a2{\alpha^2}
\def\bb{\beta}

\def\ss{\sigma}

\def\rr{\rho}

\def\pp{\partial}
\def\RN{\R^N}
\def\C{\field{C}}

\def\gd{\nabla}
\def\lp{\triangle}

\def\mb{\mathcal{B}}
\def\mc{\mathcal{C}}
\def\mh{\mathcal{H}}

\begin{document}
\title[Hamiltonian Regularization of the Nonlinear Schr\"{o}dinger Equation]
{Nonlinear Schr\"{o}dinger-Helmholtz Equation as Numerical Regularization of the Nonlinear
Schr\"{o}dinger Equation}
\date{June 27, 2007}

\author[Y.Cao]{Yanping Cao}
\address[Y.Cao]
{Department of Mathematics\\
University of California\\
Irvine, CA 92697-3875,USA}
\email{ycao@math.uci.edu}

\author[Z. H. Musslimani]{Ziad H. Musslimani}
\address[Z. H. Musslimani]
{Department of Mathematics\\
Florida State University\\
Tallahassee FL 32306, USA}
\email{musliman@mail.math.fsu.edu}

\author[E.S. Titi]{Edriss S. Titi}
\address[E.S. Titi]
{Department of Mathematics \\
and  Department of Mechanical and  Aerospace Engineering \\
University of California \\
Irvine, CA  92697-3875, USA \\
{\bf ALSO}  \\
Department of Computer Science and Applied Mathematics \\
Weizmann Institute of Science  \\
Rehovot 76100, Israel}
\email{etiti@math.uci.edu and edriss.titi@weizmann.ac.il}

\begin{abstract}
A regularized  $\alpha-$system of the Nonlinear Schr\"{o}dinger
Equation (NLS) with $2\ss$ nonlinear power in dimension $N$ is
studied. We prove existence and uniqueness of local solution in the
case $1\le \ss <\frac{4}{N-2}$ and existence and uniqueness of
global solution in the case $1\le \ss<\frac{4}{N}$. When
$\alpha\rightarrow 0^+$, this regularized system will converge to
the classical NLS in the appropriate range. In particular, the
purpose of this numerical regularization is to shed light on the
profile of the blow up solutions of the original Nonlinear
Schr\"{o}dinger Equation in the range $\frac{2}{N}\le \ss
<\frac{4}{N}$, and in particular for the critical case
$\ss=\frac{2}{N}$.
\end{abstract}

\maketitle


\section{Introduction}   \label{SEC-intro}
The Nonlinear Schr\"{o}dinger equation (NLS):
\begin{eqnarray}
&&\hskip-.8in
iv_t+ \lp v+|v|^{2\ss}v=0,\qquad \qquad x\in \RN, \quad t\in \R,    \label{NLS}\\
&&\hskip-.8in
 v(0)=v_0,            \nonumber
 \end{eqnarray}
where $v$ is a complex-valued function in $\RN\times \R$,
arises in various physical contexts describing wave
propagation in nonlinear media (see, e.g., \cite{Kelley}, \cite{Sulem}, \cite{Taniuti}
and \cite{Zakharov}). For example,
when $\sigma=1$, equation (\ref{NLS}) describes propagation of a laser beam
in a nonlinear optical medium whose index of refraction is proportional to
the wave intensity.
Also, the Nonlinear Schr\"{o}dinger Equation successfully models other wave phenomena such as
water waves at the free surface of an ideal fluid as well as plasma waves. In all cases, it is interesting to note that Eq.~(\ref{NLS}) describes wave propagation in nonhomogeneous linear media with self-induced potential  given by $|v|^{2\ss}$.\\

As it is mentioned above, the $\sigma =1$ case is particularly interesting for laser beam propagation in optical Kerr media.
Depending on the dimensionality of the space upon which the beam is propagating in, the wave dynamics can be either ``simple" or ``intricate".
In one space dimension, the NLS equation is known to be integrable and possesses soliton solution that preserves their structure upon collision
\cite{as}. The picture in two-dimensional (2D) space is totally different. The 2D NLS equation is not integrable, hence no exact soliton solutions are known. Instead, the 2D NLS equation admits the waveguide solution (also known as Townes soliton) $v(x,y,t)=R(r)\exp(it)$ with $r=\sqrt{x^2+y^2}$ where $R>0$ satisfies the nonlinear boundary value problem
\begin{equation}
\frac{d^2R}{dr^2}+\frac{1}{r}\frac{dR}{dr} -R+R^3=0\;,\;\;\;\;\; \frac{dR}{dr}(0)=0\;,\;\;\;\;
\lim _{r\rightarrow +\infty}R(r)=0\;.
\end{equation}
Importantly, the $L^2$ norm (or power in optics) of the Townes soliton defines a critical value
 for blow up. If initially the beam's power is larger than that of the Townes soliton,
$||v_0||^2_{L^2} > ||R||^2_{L^2},$ then the beam undergoes a finite time blow up. If on the
other hand
 $||v_0||^2_{L^2} < ||R||^2_{L^2}$, then the wave will diffract.
Various mechanisms to arrest collapse have been suggested such as
nonparaxiality \cite{fibich1}, or higher order dispersion \cite{fibich2}.

As a result, an important issue that arises in the mathematical
study of the NLS is the question of local and global existence of
solutions, their uniqueness, as well as the profile of blow up
solutions. Knowing answers to such questions may have some
consequences on possible physical observations of phenomenon
governed by the NLS and in validating its derivation. In the works
of Ginibre and Velo \cite{Velo} and Weinstein \cite{Weinstein}, it
is proved that equation (\ref{NLS}) has a unique global solution
when $0<\ss<\frac{2}{N}$, and that it has a unique global solution
for ``small" initial data for the critical case $\ss=\frac{2}{N}$.
The proof of global existence uses the fact that the energy
$\mathcal{N}(v)=\int_{\RN} |v(x,t)|^2\: dx$ and the Hamiltonian
$\mathcal{H}(v)=\int_{\RN} \left( | \gd v(x,t)|^2-\frac{
|v(x,t)|^{2\ss+2}}{\ss+1} \right)\: dx$ are conserved quantities of
the dynamics of (\ref{NLS}). In the case of $\ss\ge\frac{2}{N}$,
Glassey \cite{Glassey} proved that there exist solutions that
develop singularities in finite time. In recent years there was an
intensive  remarkable computational work  concerning the blow up for
the critical case $\ss=\frac{2}{N}$. For instance,  Merle and
Raphael have obtained a sharp lower bound on the blow up rate for
the $L^2$ norm of the  NLS in $\RN$ (see \cite{MerleRaphael} and
references therein). Moreover, Fibich and Merle \cite{FibichMerle}
studied self-focusing in bounded domains using a combination of
rigorous, asymptotic and numerical results.
\\\\
Instead of the potential $|v|^{2\ss}$, physicists consider self-gravational
potential (see, e.g., \cite{Penrose}, \cite{Ruffini}) and come to a new system:
Schr\"{o}dinger-Newton equation (SN):
\begin{eqnarray}
&&\hskip-.8in
iv_t+\lp v+ \psi v=0,\qquad \qquad x \in \RN, \quad t \in \R,          \label{SN}\\
&&\hskip-.8in
-\aa ^2\lp \psi = |v|^2,                                               \nonumber \\
&&\hskip-.8in
v(0)=v_0,                                                                 \nonumber
\end{eqnarray}
where $\aa>0$ is a real constant. System (\ref{SN}) is a Hamiltonian system with
the corresponding Hamiltonian
$\mathcal{H}(v)=\int_{\RN} \left( |\gd v(x,t)|^2-\frac{\psi(x,t) |v(x,t)|^2}{2} \right)\:dx$
and can be obtained formally by the variational principle
$i\frac{\pp v}{\pp t}=\frac{\delta \mathcal{H}(v)}{\delta v^*}$, where $v^*$ denotes the complex
conjugate of $v$.
System (\ref{SN}), or at least its stationary state, has been studied
\cite{MorozTod}, \cite{TodMoroz}, \cite{Lie}.
This coupled system of equations consists of the Schr\"{o}dinger equation for a wave function $v$ moving
in a potential $\psi$, where $\psi$ is obtained by solving the Poisson
equation with source $\rho=|v|^2$. It can be thought of as the Schr\"{o}dinger
equation for a particle moving in its own gravatational field.
As in the NLS, the energy $\mathcal{N}(v)=\int_{\RN} |v(x,t)|^2 \:dx $ and the Hamiltonian
$\mathcal{H}(v)=\int_{\RN} \left( |\gd v(x,t)|^2 -\frac{ \psi(x,t) |v(x,t)|^2}{2} \right) \:dx$
are also conserved in this system. The question of existence and uniqueness of local and global
solutions for system ({\ref{SN}) has not been answered completely yet.\\\\
Inspired by the $\alpha-$models of turbulence (see, e.g.,
\cite{Bardina}, \cite{NSa}, \cite{Leraya}, \cite{FHT}, \cite{HMR}, \cite{MLa} and references therein), we introduce a generalization of
(\ref{SN}),
the Schr\"{o}dinger-Helmholtz (SH) regularization
of the classical NLS:
\begin{eqnarray}
&&\hskip-.8in
iv_t+ \lp v+ u|v|^{\ss-1}v=0,\qquad \qquad x \in \RN, \quad t \in \R,   \label{SNG}\\
&&\hskip-.8in
u-\a2 \lp u=|v|^{\ss+1},                                                    \nonumber \\
&&\hskip-.8in
v(0)=v_0,                                                                  \nonumber
\end{eqnarray}
where $\aa>0$ and $\sigma \ge 1$.
System (\ref{SNG}) is a Hamiltonian system with the corresponding Hamiltonian
$\mathcal{H}(v)=\int_{\RN} \left( |\gd v(x,t)|^2-\frac{u(x,t)|v(x,t)|^{\ss+1}}{\ss+1} \right ) \: dx$
and can be obtained formally by the variational principle
$i \frac{\pp v}{\pp t}=\frac{\delta \mathcal{H}(v)}{\delta v^*}$, where again $v^*$ denotes
the complex conjugate of $v$.
In this system, we can regard the wave
function $v$ moves in a potential $u|v|^{\ss-1}$, where $u$ is obtained by
solving the Helmholtz elliptic problem $ u-\aa^2 \lp u=|v|^{\ss+1}$.
Observe that the energy
$\mathcal{N}(v)=\int_{\RN} |v(x,t)|^2 \:dx$ and the Hamiltonian
$\mathcal{H}(v)=\int_{\RN} \left( |\gd v(x,t)|^2 -\frac{u(x,t)|v(x,t)|^{\ss+1}}{\ss+1} \right) \: dx$
are conserved in this system.
When $\ss=1$, we have the potential $u$ as in the SN with the only difference that
the Possion equation is modified as a Helmholtz equation. So we consider this system as
a generalized system of SN.
A more important fact is that when $\aa=0$, one recovers the classical NLS,
therefore we regard this system as a regularization of the classical NLS.
In this paper we focus on the case $\aa>0$ and in our subsequential works, we will
investigate the behavior when $\aa \rightarrow 0^+$.
In particular, we will investigate the case $\ss=\frac{2}{N}$,
which is not completely understood.\\\\
In this paper, we will study the question of local and global existence of unique
solution for system (\ref{SNG}).
Specifically, we will prove the short time existence of unique solution,
when $1\le \ss< \frac{4}{N-2}$
(we define once and for all
$\frac{4}{N-2}=\infty$ when $N\le2$).
Moreover, we will show global existence of unique solution when
$1\le \sigma < \frac{4}{N}$. The proof will follow the ideas of \cite{Velo} and
\cite{Weinstein} and use the important fact of the conservation
of the corresponding energy and the Hamiltonian of (\ref{SNG}).
All the proofs presented here will apply directly to system (\ref{SN}) as well.
So simultaneously we have the same results for system (\ref{SN}):
we have short time existence of unique solution when
$1 \le \ss<\frac{4}{N-2}$, and global existence of unique solution
when $1\le \ss <\frac{4}{N}$. Comparing to the results of the classical
NLS (\ref{NLS}) ($\ss<\frac{2}{N-2}$ for local existence and $\ss<\frac{2}{N}$ for global existence), one
expects these ``better" results for (\ref{SN}) and (\ref{SNG})
since the nonlinear terms in (\ref{SN}) and (\ref{SNG})
are milder than that of the classical NLS (\ref{NLS}). The parametre $\aa$ plays an important role in our proofs.
In a subsequential paper,
we will investigate numerically the blow up profiles of the NLS,
in the relevant range of $\ss$, when $\aa \rightarrow 0^+$.
\\\\
In section \ref{SEC-pre}, we will introduce some essential notations and definitions, and some
preliminary results that will be used throughout the paper.
Following the work of Ginibre and Velo \cite{Velo}, we prove in section \ref{SEC-local}
local (in time) existence and uniqueness of solution for system (\ref{SNG}) using the contraction
mapping principle.
In section \ref{SEC-global}, we will extend the local solution to global existence,
for $1 \le \ss <\frac{4}{N}$, after
establishing the required \textit{a priori} estimates for the $H^1$ norm of the solution, which
remains finite for every finite interval of time.
\section{Notations and Preliminaries}   \label{SEC-pre}
In this section we introduce some preliminary results and the basic notations
and definitions that will be used
throughout this paper.
\\\\
We denote by $\| \cdot \|_p$ the norm in the space $L^p=L^p(\RN)$ ($1 \le p\le \infty)$, except
for $p=2$ where the subscript 2 will be omitted.
We will denote by $\langle \cdot,\cdot \rangle$ the scalar product in $L^2$.
The conjugate pair $p,\: p'$ satisfies the relation $\frac{1}{p}+\frac{1}{p'}=1$.
For any real number $l$, we denote by $H^l=H^l(\RN)$, the usual Sobolev space.
Of special interest is the $H^1$ Sobolev space with the norm defined by
\begin{equation}
\hskip-.8in
\|v\|_{H^1}^2=\int_{\RN} \left(1+|\xi|^2 \right) |\hat{v}(\xi)|^2 d \xi,           \label{h1-fourier}
\end{equation}
or equivelently,
\begin{equation}
\hskip-.8in
\|v\|_{H^1}^2=\|v \|^2 +\|\gd v\|^2.     \label{h1-norm}
\end{equation}
We denote by
$\|u\|_{W^{k,p}}=\left(\Sigma_{|\aa|\le k} \int_{\RN} |D^\aa u|^p \: dx \right)^{1/p},1\le p < \infty$,
for $u$ belongs to the Sobolev space $W^{k,p}(\RN)$. For any interval $I$ of the real line $\R$,
and for any Banach space $\mb$,
we denote by $\mc(I,\mb)$ (respectively $\mc_b(I, \mb)$) the space of continuous
(repectively bounded continuous) functions from $I$ into $\mb$.\\\\
In this paper $C$ and $C_\aa$ will denote constants which might depend on various parameters
of the problem. They might vary in value from one time to another, but they are independent of
the solution. When it is relevant we will comment on the asymptotic behavior of these constants
as they depend on the corresponding parameters.\\\\
First, we recall some classical
Gagliardo-Nirenberg and Sobolev inequalities (see, e.g., \cite{Adams}).
\begin{proposition}
(1) For any $N\ge 1$, we have
\begin{eqnarray}
&&\hskip-.8in
\|v\|_q \le C \|v\|^{1-\frac{q-2}{2q}N} \|\gd v\|^{\frac{q-2}{2q}N}
\qquad   \mbox{ for every } v \in H^1,\: 0<\frac{q-2}{2q}N\le 1    \label{H-1}  \\
&&\hskip-.8in
\| v \| _q \le C\| v \| _{W^{2,m}}  \qquad \qquad \qquad \qquad
\mbox{for every } v \in W^{2,m},\: q \ge m,\: 2m >N                    \label{2-m-m} \\
&&\hskip-.8in
\|v\|_q \le C\|v\|_{W^{2,m}}  \qquad \qquad \qquad \qquad
\mbox{for every } v \in W^{2,m},\: \frac{1}{q} \ge \frac{1}{m}-\frac{2}{N}\ge 0,\: q<\infty  \label{w-2-p}\\
\mbox{In particular},\\
&&\hskip-.8in
\|v\|_q \le C\|v\|_{W^{2,2}}=C\|v\|_{H^2} \qquad \quad \quad
\mbox{ for every } v \in H^2,\: 2\le q\le \infty,\: N\le3.                       \label{2-2-q}
\end{eqnarray}
(2) For $N\le 2$,
\begin{equation}
\hskip-.8in
\|v\|_q \le C \|v\|_{H^1} \qquad \qquad \qquad \qquad
\mbox{ for every } v \in H^1, \:2\le q< \infty.                              \label{H-1-q}   \\
\end{equation}
\end{proposition}

With these inequalities at hand, we can process the nonlinear term.
Let us rewrite the term
\begin{equation}
f(v)=u|v|^{\ss-1}v=B(|v|^{\ss+1})|v|^{\ss-1}v,    \label{fv}
\end{equation}
where $B=(I-\aa^2 \lp)^{-1}$, the inverse of the Helmholtz operator.
Then $f$ is a locally Lipschitz mapping  from $H^1$ into $L^{r'}$, for
some $ r  \in (2, \frac{2N}{N-2}]$, where $\frac{1}{r}+\frac{1}{r'}=1$.
\begin{proposition}
Let $N\ge1$ and $1\le \ss <\frac{4}{N-2}$. For every $ v_1, v_2 \in H^1 \subset L^r$,
where $r$ depends on the given $\ss$ and belongs to the range $r \in (2,\frac{2N}{N-2} ]$
(we consider $\frac{2N}{N-2}$ as $\infty$ when $N\le2$), we have
$\|f(v_1)-f(v_2)\|_{r'}\le k \|v_1-v_2\|_r$,
where $k=C_\aa \left( \|v_1\|_{H^1}+\|v_2\|_{H^1} \right)^{2\ss}$ and
$\frac{1}{r}+\frac{1}{r'}=1$, for some constant $C_\aa$.
\end{proposition}
Before we prove this proposition, we will state the following Lemmas:
\begin{lemma}
Let $N\ge 1$ and $1\le \ss < \frac{4}{N-2}$. For every
$ v_1, v_2, v \in H^1 \subset L^r$, where $r$ depends on the given $\ss$
and belongs to the range $r \in (2,\frac{2N}{N-2}]$, we have
\begin{eqnarray}
&&\hskip-.8in
\|B( |v|^\ss) |v|^{\ss-1} v (v_1-v_2)\|_{r'}
\le C_\aa\ \|v\|_{H^1}^{2\ss}\|v_1-v_2\|_r\:,             \label{I-1} \\
&&\hskip-.8in
\|B(|v|^{\ss-1}v)|v|^{\ss-1}v(v_1-v_2)\|_{r'}
\le C_\aa \|v\|_{H^1}^{2\ss} \|v_1-v_2\|_r\:.             \label{I-1-2}
\end{eqnarray}
\end{lemma}
\begin{proof}
First, denote
\begin{equation*}
\hskip-.8in
I_1=\|B(|v|^{\ss})|v|^{\ss-1}v(v_1-v_2)\|_{r'}\:.
\end{equation*}
{\bf{Case 1. $N\le 2$:}}\\
By H\"{o}lder's inequality, we have
\begin{eqnarray}
\hskip-.8in
I_1
&\le & \|B ( |v|^{\ss})|v|^{\ss-1}v\|_{r'\bb_1} \|v_1-v_2\|_{r'\gamma_1}   \nonumber  \\
&=& \|B(|v|^\ss) |v|^{\ss-1} v\|_{\frac{r}{r-2}} \|v_1-v_2\|_r\:,          \label{lemma1-1}
\end{eqnarray}
where in the last equality, we choose $\gamma_1=r-1>1$ and
$\frac{1}{\bb_1}+\frac{1}{\gamma_1}=1$ such that
$r'\gamma_1=r$ and $r'\bb_1=\frac{r}{r-2}$.
\\\\
By Cauchy-Schwarz inequality, we have
\begin{equation}
\hskip-.8in
\| B( |v|^{\ss} ) |v|^{\ss-1}v\|_{\frac{r}{r-2}}
\le \|B(|v|^{\ss})\|_{\frac{2r}{r-2}} \| |v|^{\ss} \|_{\frac{2r}{r-2}}\:.   \label{lemma1-2}
\end{equation}
Now, for the elliptic equation $u-\aa ^2 \lp u=f $ in $\RN$, we have the
regularity property \cite{func}, \cite{Yudovich}
\begin{equation}
\hskip-.8in
\|u\|_{W^{2,p}} \le C_\aa \|f\|_p    \qquad \qquad  \mbox{for any } 1< p<\infty,  \label{elliptic}
\end{equation}
where $C_\aa$ depends on $N,p$ and $\aa$, and $C_\aa \sim \frac{1}{\aa ^2}$
as $\aa \rightarrow 0^+$. Moreover, for $\aa$ fixed, $C_\aa \sim p $ as $p \rightarrow \infty$.\\
Since $\frac{2r}{r-2}>2$, by (\ref{2-2-q}) and (\ref{elliptic}), we have
\begin{eqnarray*}
\| B(|v|^{\ss})\|_{\frac{2r}{r-2}}
&\le& C \|B (|v|^{\ss})\|_{W^{2,2}}  \\
&\le& C_\aa  \| |v|^{\ss}\|               \\
&=& C_\aa  \|v\|^{\ss}_{2\ss}\:,
\end{eqnarray*}
and
\begin{equation*}
\hskip-.8in
\| |v|^{\ss} \|_{\frac{2r}{r-2}}=\|v\|_{\frac{2r}{r-2}\ss}^{\ss}\:.
\end{equation*}
Since $\ss\ge 1$, combining the above two terms
and applying (\ref{H-1-q}), we have
\begin{equation*}
\hskip-.8in
I_1\le C_\aa  \|v\|_{H^1}^{2\ss} \|v_1-v_2\|_r\:.
\end{equation*}
\vskip.1in
\noindent
{\bf{Case 2. $N\ge3$:}}\\
Applying H\"{o}lder's inequality, we have
\begin{eqnarray*}
\hskip-.8in
I_1
&\le& \| B(|v|^{\ss})\|_{r'\theta_1} \| |v|^{\ss-1}v\|_{r'\bb_1} \| v_1-v_2\|_{r'\gamma_1}  \\
&=& \|B (|v|^{\ss})\|_{r'\theta_1} \|v\|_{\ss r'\bb_1}^{\ss} \|v_1-v_2\|_{r'\gamma_1}\:,
\end{eqnarray*}
where $\frac{1}{\theta_1}+\frac{1}{\bb_1}+\frac{1}{\gamma_1}=1$.\\
Now by (\ref{2-m-m}), (\ref{w-2-p}) and (\ref{elliptic}), we have
\begin{eqnarray*}
\hskip-.8in
\| B(|v|^{\ss})\|_{r'\theta_1}
&\le& C \| B(|v|^{\ss})\|_{W^{2,m}}   \\
&\le& C_\aa \| |v|^{\ss}\|_m       \\
&=& C_\aa \|v\|_{\ss m}^{\ss}\:,
\end{eqnarray*}
where we require  $\frac{1}{r'\theta_1}\ge \frac{1}{m}-\frac{2}{N}$ when
$\frac{1}{m}-\frac{2}{N}\ge 0$, or
$r'\theta_1 \ge m $ when $\frac{1}{m}-\frac{2}{N} <0$, for $m>1$ to be determined later.\\

Therefore, we obtain
\begin{equation}
\hskip-.8in
I_1 \le C_\aa \|v\|_{\ss m}^{\ss} \|v\|_{\ss r'\bb_1}^{\ss} \|v_1-v_2\|_{r'\gamma_1}\:. \label{I-1-simple}
\end{equation}
Now, by requiring $\ss m=\ss r'\bb_1=r'\gamma_1=r$, we have
\begin{eqnarray*}
&&\hskip-1in
\theta_1=\frac{r-1}{r-\ss-2}>1 \Rightarrow \ss<r-2    \\
&&\hskip-1in
\bb_1=\frac{r-1}{\ss} >1   \Rightarrow \ss<r-1   \\
&&\hskip-1in
\gamma_1=r-1>1                                    \\
&&\hskip-1in
m=\frac{r}{\ss}>1                                    \\
&&\hskip-1in
\ss <\frac{N+2}{2N}r-1 \quad \mbox{or } \quad \ss >\frac{r-2}{2}\:.
\end{eqnarray*}
Since $2 < r \le \frac{2N}{N-2}$, we conclude that $\ss < \frac{4}{N-2}$, i.e.,
\begin{equation*}
\hskip-.8in
I_1 \le C_\aa \|v\|_{H^1}^{2\ss} \|v_1-v_2\|_r\:.
\end{equation*}
By exactly the same steps, inequality (\ref{I-1-2}) follows readily.
\end{proof}
\begin{lemma}
Let $N\ge 1$ and $1\le \ss <\frac{4}{N-2}$\:. For every $v_1,v_2,v \in H^1 \subset L^r$,
where $r$ depends on the given $\ss$ and belongs to the range $ r\in (2,\frac{2N}{N-2}]$,
we have
\begin{eqnarray}
&&\hskip-.8in
\|B (|v|^{\ss+1})|v|^{\ss-1}(v_1-v_2)\|_{r'}
\le C_\aa \|v\|^{2\ss}_{H^1} \|v_1-v_2\|_r\:,          \label{I-2}   \\
&&\hskip-.8in
\| B(|v|^{\ss+1}) |v|^{\ss-3}v^2(v_1-v_2) \|_{r'}
\le C_\aa \|v\|^{2\ss}_{H^1} \|v_1-v_2\|_r\:.          \label{I-2-2}
\end{eqnarray}
\end{lemma}
\begin{proof}
Denote
\begin{equation*}
\hskip-.8in
I_2=\| B(|v|^{\ss+1}) |v|^{\ss-1} (v_1-v_2)\|_{r'}\:.
\end{equation*}
{\bf{Case 1. $ N \le 2$:}}\\
By H\"{o}lder's inequality, we obtain
\begin{eqnarray*}
I_2
&\le& \| B(|v|^{\ss+1}) |v|^{\ss-1} \|_{r'\bb_1} \|v_1-v_2\|_{r'\gamma_1}  \\
&=& \| B(|v|^{\ss+1})|v|^{\ss-1}\|_{\frac{r}{r-2}}\|v_1-v_2\|_r\:,
\end{eqnarray*}
where we choose the same $\bb_1,\gamma_1$ as in (\ref{lemma1-1}).\\\\
Now, when $\ss >1$, again by H\"{o}lder's inequality, we have
\begin{eqnarray}
\| B(|v|^{\ss+1})|v|^{\ss-1} \|_{\frac{r}{r-2}}
&\le& \| B(|v|^{\ss+1})\|_{\frac{r}{r-2}\bb_2} \| |v|^{\ss-1}\|_{\frac{r}{r-2}\gamma_2}  \nonumber \\        \\
& = & \| B(|v|^{\ss+1})\|_{\frac{r}{r-2}\bb_2} \|v\|_{(\ss -1)\frac{r}{r-2}\gamma_2}^{\ss-1}\:. \label{new}
\end{eqnarray}
By choosing $1<\bb_2<\ss$ and $2<r \le 4$, one can easily verify that
$\frac{r}{r-2}\bb_2 \ge 2$ and $(\ss-1)\frac{r}{r-2}\gamma_2 >2$.
\\\\By (\ref{2-2-q}) and (\ref{elliptic}), we obtain
\begin{eqnarray}
\| B(|v|^{\ss+1})\|_{\frac{r}{r-2}\bb_2}
&\le& C \| B(|v|^{\ss+1})\|_{W^{2,2}}                 \nonumber   \\
&\le& C_\aa \| |v|^{\ss+1} \|                         \nonumber   \\
&=& C_\aa \| v\|_{2(\ss+1)}^{\ss+1}\:.                   \label{new2}
\end{eqnarray}
Since $2(\ss+1)>2$ and $(\ss-1)\frac{r}{r-2}\gamma_2>2$ when $ \ss >1$, by (\ref{H-1-q})
we conclude
\begin{equation}
I_2 \le C_\aa \|v\|_{H^1}^{2\ss}\|v_1-v_2\|_r\:.
\end{equation}
Now, when $\ss=1$, by choosing $2<r\le4$ in (\ref{new}) and applying
inequality (\ref{new2}) with $\beta_2=1$,
we have
\begin{eqnarray*}
\|B(|v|^{\ss+1} ) |v|^{\ss-1}\|_{\frac{r}{r-2}}
&=&  \| B(|v|^{\ss+1} )\|_{\frac{r}{r-2}}   \\
&\le & C_{\alpha} \| v\|_{2(\ss+1)}^{\ss+1}\:.
\end{eqnarray*}
By (\ref{H-1-q}), we conclude that
\begin{equation}
I_2 \le C_\aa \|v\|_{H^1}^{2\ss}\|v_1-v_2\|_r\:.
\end{equation}
\vskip.1in
\noindent
{\bf{Case 2. $N\ge 3$:}}\\
By H\"{o}lder's inequality, we have
\begin{eqnarray}
I_2
&\le& \| B(|v|^{\ss+1})\|_{r'\theta_1} \| |v|^{\ss-1} \|_{r'\bb_1} \|v_1-v_2\|_{r'\gamma_1}  \nonumber \\
&=&\| B(|v|^{\ss+1})\|_{r'\theta_1}\|v\|_{(\ss -1)r'\bb_1}^{\ss-1} \|v_1-v_2\|_{r'\gamma_1}\:,  \label{Holder}
\end{eqnarray}
where $\frac{1}{\theta_1}+\frac{1}{\bb_1}+\frac{1}{\gamma_1}=1$.\\

Now, by (\ref{2-m-m}), (\ref{w-2-p}) and (\ref{elliptic}), we have
\begin{eqnarray*}
\hskip-.8in
\| B(|v|^{\ss+1})\|_{r'\theta_1}
&\le& C\| B(|v|^{\ss+1})\|_{W^{2,m}}                            \\
&\le& C_\aa \| |v|^{\ss+1} \| _m                               \\
&=& C_\aa \| v\| _{(\ss+1)m}^{\ss+1}\:,
\end{eqnarray*}
where $\frac{1}{r'\theta_1}\ge \frac{1}{m}-\frac{2}{N}$ when $\frac{1}{m}-\frac{2}{N} \ge 0$,
and $r'\theta_1 \ge m$ when $\frac{1}{m}-\frac{2}{N}<0$, and $m>1$ to be decided later.\\
Then we obtain
\begin{equation}
\hskip-.8in
I_2\le C_\aa\|v\|_{(\ss+1)m}^{\ss+1} \|v\|_{(\ss-1)r'\bb_1}^{\ss-1} \|v_1-v_2\|_{r'\gamma_1}\:.  \label{I-2-simple}
\end{equation}
Now, by requiring $(\ss+1)m=(\ss-1)r'\bb_1=r'\gamma_1=r$
 (choose $\bb_1=\infty$ when $\ss=1$), we have
\begin{eqnarray*}
&&\hskip-1in
\theta_1=\frac{r-1}{r-1-\ss}>1 \Rightarrow \ss<r-1  \\
&&\hskip-1in
\bb_1=\frac{r-1}{\ss-1}>1   \Rightarrow \ss<r    \\
&&\hskip-1in
\gamma_1=r-1>1               \\
&&\hskip-1in
m=\frac{r}{\ss+1}>1            \\
&& \hskip-1in
\ss<\frac{N+2}{2N}r-1 \quad \mbox{or }\quad  \ss>\frac{r-1}{2}\:.
\end{eqnarray*}
Since $2< r \le \frac{2N}{N-2}$, we obtain $\ss<\frac{4}{N-2}$, i.e.,
\begin{equation*}
\hskip-.8in
I_2 \le C_\aa \|v\|_{H^1}^{2\ss}\|v_1-v_2\|_r\:.
\end{equation*}
The same can be shown for (\ref{I-2-2}).
\end{proof}
\vskip.1in
Now, we are ready to prove Proposition 2.

\begin{proof}
Recall that $f(v)=B(|v|^{\ss+1})|v|^{\ss-1}v$, by direct calculation, we have
\begin{eqnarray}
&&\hskip.8in
\frac{\partial f(v)}{\partial v}=
\frac{\ss+1}{2}[B(|v|^{\ss})|v|^{\ss-1}v+B(|v|^{\ss+1})|v|^{\ss-1} ]  \label{pfpv} \\
&&\hskip.8in
\frac{\partial f(v)}{\partial v^*}=
\frac{\ss+1}{2} B(|v|^{\ss-1}v)|v|^{\ss-1}v+\frac{\ss-1}{2}B(|v|^{\ss+1})|v|^{\ss-3}v^2 \label{pfpvstar}
\end{eqnarray}
when $v\ne0$, where $v^*$ is the conjugate of the complex-valued function $v$.\\
For $v=0$, we have
\begin{eqnarray*}
&&\hskip-.8in
\frac{\partial f(v)}{\partial v}=0,                 \\
&&\hskip-.8in
\frac{\partial f(v)}{\partial v^*}=0.
\end{eqnarray*}
Now, by Lemma 3 and Lemma 4 we obtain
\begin{eqnarray*}
&&\hskip-.8in
\|\frac{\partial f(v)}{\partial v}(v_1-v_2)\|_{r'} \le C_\aa \|v\|_{H^1}^{2\ss} \|v_1-v_2\|_r\:, \\
&&\hskip-.8in
\|\frac{\partial f(v)}{\partial v^*} (v_1-v_2)\|_{r'} \le C_\aa \|v\|_{H^1}^{2\ss}\|v_1-v_2\|_r\:.
\end{eqnarray*}

When $v_1v_2\ne 0$, by Mean-Value Theorem, we have
\begin{equation*}
|f(v_1)-f(v_2)| \le
\max \{  |\frac{\partial f(v)}{\partial v}|,|\frac{\partial f(v)}{\partial v^*}| \}
 |v_1-v_2|.
\end{equation*}
That is, for some intermediate point $v_0$ between $v_1$ and $v_2$
\begin{eqnarray*}
\|f(v_1)-f(v_2)\|_r
&\le& {\tilde{C_\aa}} \|v_0\|_{H^1}^{2\ss} \|v_1-v_2\|_r   \\
&\le& C_\aa \left( \|v_1\|_{H^1}+\|v_2\|_{H^1} \right) ^{2\ss} \|v_1-v_2\|_r\:.
\end{eqnarray*}

When $v_1v_2=0$, without loss of generality, we assume that $v_2=0$, then
\begin{eqnarray*}
\|f(v_1)-f(v_2)\|_{r'}
&=&\|f(v_1)-0\|_{r'}            \\
&=&\| B(|v_1|^{\ss+1})|v_1|^{\ss-1}v_1\|_{r'}  \\
&\le& C_\aa \left( \|v_1\|_{H^1}+\|v_2\|_{H^1} \right)^{2\ss} \|v_1\|_r\:,
\end{eqnarray*}
where the last inequality follows from direct application of Lemma 4.\\

Therefore, we conclude that, for any $v_1,v_2 \in H^1 \subset L^r$,
\begin{equation*}
\hskip-.8in
\|f(v_1)-f(v_2)\|_{r'} \le C_\aa
\left( \|v_1\|_{H^1}+\|v_2\|_{H^1} \right)^{2\ss} \|v_1-v_2\|_r\:.
\end{equation*}
\end{proof}

Next, we will give some elementary properties of the free evolution (linear
Schr\"{o}dinger equation) formally
defined by the group of operators
\begin{equation}
\hskip-.8in
U(t)=\exp (it \lp) \label{U},
\end{equation}
where $t \in \R$. In the following, we will
state some well-known results about the operator $U(t)$ without proving them (see, e.g.,
\cite{Velo},\cite{Sulem}).\\
\begin{lemma}
For any $r\ge 2$,  and for any $t\neq 0$,
$U(t)$ is a bounded linear operator from $L^{r'}$ to $L^r$, and the map $t \rightarrow U(t)$
is strongly
continous. Moreover, for all $t \in \R\setminus\{0\}$, one has
\begin{equation}
\hskip-.8in
\|U(t)v\|_r \le (4\pi |t|)^{\frac{N}{r}-\frac{N}{2}} \|v\|_{r'}  \label{U-conts}
\end{equation}
for all $ v \in L^{r'}$.
\end{lemma}
\begin{corollary}
Let $I$ be an interval of $\R$, and let $v \in \C(I,L^{r'})$. Then for all $t \in \R$ the map
$\tau \rightarrow U(t-\tau)v(\tau)$ is continous from $I\setminus \{t\}$ into $L^r$.
\end{corollary}

\section{Existence and Uniqueness of Local Solutions} \label{SEC-local}
In this section we will prove a local existence and uniquenss theorem
of solutions to system (\ref{SNG}) by a fixed point technique.\\
The integral equation
\begin{equation}
\hskip-.8in
v(t)=U(t-t_0)v_0+i\int_{t_0}^t U(t-\tau ) f(v(\tau )) \:d\tau  \label{integral}
\end{equation}
may be considered as the integral version of the initial value problem for
equation (\ref{SNG}).
\\\\
Defining the subspace $Y(I) \subset \mc(I,X)$ and $Y_b(I) \subset \mc_b(I,X)$ by
\begin{eqnarray*}
&&\hskip-.8in
Y(I)=\{ v:v \in \mc (I,X) \quad \mbox{and } v(t)=U(t-s)v(s)\quad \mbox{for all $s$ and $t$}\in I \}   \\
&&\hskip-.8in
Y_b(I)=Y(I) \cap\mc _b(I,X).
\end{eqnarray*}
Here for special interest we choose the Banach space $X=L^r(\RN)$, for some $r>2$, which is
specified in the proof of Lemma 3 and Lemma 4,
and $\bar{X}=L^{r'}(\RN)$.

If $v \in \mc_b(I,X)$, we shall denote its norm by $|v|_I$, and for $v \in \mc_b(I,H^1)$,
we denote its norm by $|v|_{H^1,I}$. The ball
of radius $R$ in $\mc_b(I,X)$ will be denoted by $B(I,R)$. \\\\

Let $t_1,t_2 \in \R$ and let $v(t)$ be a family of complex-valued
functions defined on $\RN$, depending on a parameter $t \in \R$. We formally
define the operators
\begin{equation}
[G(t_1,t_2)v](t)=i\int_{t_1}^{t_2}U(t-\tau)f(v(\tau))\: d\tau,   \label{G-operator}
\end{equation}
where $f$ is the nonlinear term defined in (\ref{fv}). The first lemma below gives a meaning to
the expression defined by (\ref{G-operator}) and contains some of its properties.
\begin{lemma}
For any interval $I\subset \R$ (possibly unbounded), the maps $(t_1,t_2,v) \rightarrow G(t_1,t_2)v$ are continuous from
$I \times I \times \mc(I,X)$ to $Y_b(\R)$. Moreover, for any $t_1,t_2 \in I,(t_1<t_2)$,
for any compact sub-interval $J$ such that $[t_1,t_2]\subset J\subset I$, and for any $t \in
[t_1,t_2]$, for any
$v_1,v_2 \in \mc(I,X)$ the G operator satisfies the estimates
\begin{equation*}
\|G(t_1,t_2)v_1(t)-G(t_1,t_2)v_2(t)\|_r \le k' \|v_1-v_2\|_J |t_2-t_1|^{\frac{N}{r}-\frac{N}{2}+1}
\end{equation*}
where $k'=k(4\pi)^{\frac{N}{r}-\frac{N}{2}},k=C_\aa \left(\|v_1\|_{H^1}+\|v_2\|_{H^1}\right)$, which is derived in the proof of Proposition 2.
\end{lemma}
\begin{proof}
For any $v \in \mc(I,X)$ the function $\tau \rightarrow f(v(\tau))$ belongs
to $\mc(I,\bar{X})$ as consequence of Proposition 2. Therefore, by Lemma 3,
for any $t \in \R\setminus \{ t \}$ the function
\begin{equation}
\hskip-.8in
\tau \rightarrow U(t-\tau)f(v(\tau))  \label{inside-integral}
\end{equation}
is continous from $I$ to $X$. To check the integrability of the function (\ref{inside-integral})
it will be enough to show the integrability of its norm.
More generally one is interested in
the integrability of
\begin{equation}
\hskip-.8in
\|U(t-\tau)[f(v_1(\tau ))-f(v_2(\tau ))]\|_r\:,   \label{inside-diff}
\end{equation}
for any $v_1,v_2 \in \mc(I,H^1) \subset \mc(I,X)$.\\

This is a direct consequence of Proposition 2 and Lemma 3: For $t \in \R$, for
every compact sub-interval $J\subset I$ and $\tau \in J$, we have
\begin{equation*}
\|U(t-\tau)[f(v_1(\tau ))-f(v_2(\tau ))]\|_{r'} \le (4\pi |t-\tau |)^{\frac{N}{r}-\frac{N}{2}} k \|v_1-v_2\|_J\:.
\end{equation*}
Finally, we come to the conclusion that
\begin{equation*}
\hskip-.8in
\|G(t_1,t_2)v_1(t)-G(t_1,t_2)v_2(t)\| \le k' \|v_1-v_2\|_J|t_2-t_1|^{\frac{N}{r}-\frac{N}{2}+1}.   \\\\
\end{equation*}
\end{proof}
Now, in order to study the equation (\ref{integral}) one needs the operators
\begin{equation}
\hskip-.8in
[F(t_0)v](t)=[G(t_0,t)v](t).   \label{F-operator}
\end{equation}
The existence and properties of $F$ follow immediately from Lemma 5.\\

For every $v \in \mc(I, H^1)\subset \mc(I,X)$,
\begin{equation}
\hskip-.8in
[A(t_0,v_0)v](t)=[F(t_0)v](t)+U(t-t_0)v_0   \label{A-defi}
\end{equation}
is a continuous map from
$\mc(I,H^1)\subset \mc(I,X)$ into $\mc(I,X)$.\\

With these notations equation (\ref{integral}) may be rewritten as
\begin{equation}
\hskip-.8in
A(t_0,v_0)v=v.    \label{A-mapping}
\end{equation}
The next lemma gives some elementary properties of the solutions of equation
(\ref{integral}). In particular, it expresses the consistency of the change of the initial time
$t_0$.
\begin{lemma}
Let I and J be two intervals of $\R$, $J\subset I$, let $t_0 \in J$, let $v_0\in H^1$ be
such that the function $t \rightarrow U(t-t_0)v_0$ belongs to $Y(I)$, and let $v \in \mc(J,X)$
be a solution of the equation (\ref{A-mapping})\\
$(\textit{i})$ The function
\begin{equation}
\hskip-.8in
\phi(v):s \rightarrow U(\cdot-s)v(s)=[\phi(v)](s)   \label{phi-operator}
\end{equation}
belongs to $\mc(J,Y(I))$ and satisfies for all $s,s'\in J$ the equality
\begin{equation}
\hskip-.8in
[\phi (v) ](s)-[\phi (v) ](s')=G(s',s)v.       \label{phi-G}
\end{equation}
Furthermore, if for some $s \in J,[ \phi (v)](s) \in Y_b(I)$, then $\phi (v)\in \mc(J,Y_b(I))$.
If in addition $J$ is bounded, then $\phi(v) \in \mc_b(J,Y_b(I))$.\\
$(\textit{ii})$ For any $s\in J$, $u$ satisfies the equation
\begin{equation}
\hskip-.8in
A(s,v(s))v=v.   \label{Av-v}
\end{equation}
\end{lemma}
\begin{proof}
Apply the operator $U(t-s)$ to equation $[A(t_0,v_0)v](s)=v(s)$ and use
the fact that
$U(t-s)[G(t_1,t_2)v](s)=[G(t_1,t_2)v](t)$ (for the proof of this identity,
we refer to Ginibre and Velo \cite{Velo})
yields
\begin{equation}
\hskip-.8in
[\phi(v) ](s)=U(\cdot -t_0)v_0+G(t_0,s)v  \label{phi-2}.
\end{equation}
>From which (\ref{phi-G}) follows immediately. The continuity properties of the
left-hand side of (\ref{phi-2}) are then a consequence of the assumptions made on
$v_0$ and of Lemma 7. Finally, putting $s'=t$ in (\ref{phi-G}) and taking the values of both members at $t$ one
obtains equation (\ref{Av-v}) at time $t$.
\end{proof}
We are now ready to discuss the problem of the existence and uniqueness of solutions
of equation (\ref{A-mapping}).
\begin{theorem}
For any $\rr>0$, there exists a  $T_0(\rr)>0$, depending only
on $\rr$, such that for any $t_0 \in \R$ and for any $v_0 \in H^1$,
for which $\|v_0\|_{H^1} \le \rr$,
equation (\ref{A-mapping}) has a unique solution on $\mc(I,X)$, where
$I=[t_0-T_0(\rr),t_0+T_0(\rr)]$ and $X=L^r$.
\end{theorem}
\begin{proof}
Let $\rr$ be a fixed postive number, let $t_0$ and $T \in \R$, $T_0>0$, and
let $I=[t_0-T,t_0+T]$. Then for every $v_1,v_2 \in H^1$ and $\|v_1-v_2\|_{H^1} \le 2\rr$,
Lemma 5 and (\ref{F-operator}) yield the inequality
\begin{equation}
\hskip-.8in
|F(t_0)v_1-F(t_0)v_2|_I \le 2k' |t-t_0|^{\frac{N}{r}-\frac{N}{2}+1}|v_1-v_2|_I \:. \label{exis-1}
\end{equation}
In particular, if we take $T=T_0(\rr)$ with $T_0(\rr)$ defined by
\begin{equation}
\hskip-.8in
4k'|T_0(\rr)-t_0|^{\frac{N}{r}-\frac{N}{2}+1}=1  \label{exis-2}
\end{equation}
in equality (\ref{exis-1}) it gives
\begin{equation}
\hskip-.8in
|F(t_0)v_1-F(t_0)v_2|_I \le \frac{1}{2}|v_1-v_2|_I \:.  \label{exis-3}
\end{equation}
Let now $v_0 \in X$ be such that $U(\cdot-t_0)v_0 \in B(I,\rr)$. Definition
(\ref{A-defi}) and estimate (\ref{exis-3}) imply
\begin{equation}
\hskip-.8in
|A(t_0,v_0)v|_I \le 2\rr \label{exis-4},
\end{equation}
and
\begin{equation}
\hskip-.8in
|A(t_0,v_0)v_1-A(t_0,v_0)v_2|_I \le \frac{1}{2}|v_1-v_2|_I\: , \label{exis-5}
\end{equation}
for all $v,v_1,v_2 \in B(I,2\rr)$, from which it follows that $A(t_0,v_0)$ is a contraction
from the ball $B(I,2\rr)$ into itself. The result is now a consequence of the contraction
mapping theorem.
\end{proof}
\section{Global Existence of Solutions}\label{SEC-global}
In this section we will study global existence of solutions
to system (\ref{SNG}) under the condition of $\ss\ge1$.
We will show below that we have global solutions
when $1\le \ss <\frac{4}{N}$. Comparing this to the results
of the classical NLS ($\ss<\frac{2}{N}$), we ``gain"
global regularity for larger range of values of $\ss$.
As we stated in the introduction, system (\ref{SNG}) will
recover the classical NLS as the parameter $\aa \rightarrow 0^{+}$,
in a subsequencial paper, we will study numerically the blow up profile
of the classical NLS by focusing on system (\ref{SNG})
with $\frac{2}{N} \le \ss< \frac{4}{N}$ when $N\le 3$.
To be more specific, the profile of blow-up in the critical
case $\ss=\frac{2}{N}$ in the classical NLS has not been known completely,
in a subsequential work, we will compute SH system (\ref{SNG}) and try to
find out the blow up profile
by forcing the parameter $\ss$ to approach zero.
\begin{theorem}
Let $v_0 \in H^1(\RN)$. If $1\le\ss <\frac{4}{N}$,  then there exists a unique solution \\
$ v \in C((-\infty,\infty);H^1(\RN))$ of the
initial-value problem (\ref{SNG}), in the sense of the equivalent integral equation.\\
Furthermore, as long as $v(x,t)$ remains in $H^1(\RN)$, the energy
\begin{equation}
\hskip-.8in
\mathcal{N}(v)=\int_{\RN}|v(x,t)|^2 \:dx   \label{energy}
\end{equation}
and Hamiltonian
\begin{equation}
\hskip-.8in
\mathcal{H}(v)=\int_{\RN} \left( |\gd v(x,t)|^2 -\frac{u(x,t)|v(x,t)|^{\ss+1}}{\ss+1} \right)\: dx
\label{hamiltonian}
\end{equation}
remain constant in time.
\end{theorem}
In the local existence theorem in section \ref{SEC-local}, we have shown that
the length $T_0$, of the interval of existence $[t_0,t_0+T_0]$, can
be taken to depend only on $\|v_0\|_{H^1}$. It follows that if $v(x,t)$ is a maximally
defined solution on $[t_0,T_{\mbox{max}})$, then either
\begin{equation*}
\hskip-.8in
T_{\mbox{max}}=+\infty
\end{equation*}
or
\begin{equation*}
\hskip-.8in
\lim _{t\rightarrow T_{\mbox{max}}^{-}} \| v(t)\|_{H^1}=+ \infty.
\end{equation*}
The heart of the global existence proof lies in the use of the invariants (\ref{energy})
and (\ref{hamiltonian}), which enable us to obtain an \textit{a priori} bound of the following type:
\begin{equation}
\hskip-.8in
\|v(x,t)\|_{H^1} \le C(\mathcal{N},\mathcal{H}).   \label{global-1}
\end{equation}
\begin{proof}
We proceed as follows:\\
>From (\ref{hamiltonian}) we have
\begin{equation}
\hskip-.8in
\|\gd v(x,t)\|^2 \le \mh +\frac{1}{\ss+1}\int_{\RN} u|v|^{\ss+1}\: dx  \label{gdv}.
\end{equation}
Observe that
\begin{eqnarray}
\hskip-.8in
\left| \int_{\RN} u|v|^{\ss+1} \:dx \right|
&\le& \|u\|_p\| |v|^{\ss+1}\|_{p'}                    \nonumber  \\
&=&\|u\|_p \|v\|_{p'(\ss+1)}^{\ss+1}                  \label{NL-term}
\end{eqnarray}
where
\begin{equation*}
\hskip-.8in
\frac{1}{p}+\frac{1}{p'}=1,\quad 1<p,\: p'<\infty.
\end{equation*}
\vskip.1in
\noindent
{\bf{Case 2. $N \le2 $:}}\\
By (\ref{2-m-m}) and (\ref{elliptic}), we have
\begin{eqnarray*}
\hskip-.8in
\|u\|_p
&\le & C \|u\|_{W^{2,p}}  \\
&=& C \| B(|v|^{\ss+1}) \|_{W^{2,p}} \\
&\le & C_\aa  \| |v|^{\ss+1}\| _p         \\
&=& C_\aa \|v\|_{(\ss+1)p}^{\ss+1}
\end{eqnarray*}
for any $p>1$.\\
Plug the above inequality into (\ref{NL-term}) and requiring $p=p'=2$, we obtain
\begin{equation*}
\hskip-.8in
\left| \int_{\RN} u|v|^{\ss+1}dx\right| \le C_\aa \|v\|_{2(\ss+1)}^{2(\ss+1)}\:.
\end{equation*}
By (\ref{H-1}), taking $q=2(\ss+1)$, we obtain
\begin{equation*}
\|v\|_{2(\ss+1)}
\le C\|v\|^{1-\frac{2(\ss +1)-2}{2\cdot2(\ss +1)}N}
\|\gd v\|^{\frac{2(\ss +1)-2}{2\cdot 2(\ss +1)}N}\:,
\end{equation*}
with
\begin{equation*}
\hskip-.8in
0<\frac{2(\ss+1)-2}{2\cdot 2(\ss+1)}N <1,
\end{equation*}
which is always satisfied when $N\le 2$.\\\\
Then (\ref{gdv}) yields
\begin{equation}
\hskip-.8in
\|\gd v(t) \|^2 \le \mathcal{H}_0 +C_\aa \|v_0\|^{2(\ss+1)-(\ss N-2)} \| \gd v\|^{\ss N-2}\:. \label{final-1}
\end{equation}
\\\\
\vskip-.1in
\noindent
{\bf{Case 2. $N \ge 3 $:}}\\
By (\ref{w-2-p}) and (\ref{elliptic}), we obtain
\begin{eqnarray*}
\hskip-.8in
\|u\|_p
&\le & C \|u\|_{W^{2,m}}  \\
&=& C \| B(|v|^{\ss+1}) \|_{W^{2,m}} \\
&\le & C_\aa  \| |v|^{\ss+1}\| _m         \\
&=& C_\aa \|v\|_{(\ss+1)m}^{\ss+1}\:,
\end{eqnarray*}
where
\begin{equation}
\hskip-.8in
\frac{1}{p}=\frac{1}{m}-\frac{2}{N}>0  \Rightarrow m<\frac{N}{2}\:.   \label{condition-2}
\end{equation}
Pluging into (\ref{NL-term}) and requiring $m=p'$, i.e., $m=\frac{2N}{N+2}$, we get
\begin{equation*}
\hskip-.8in
\left|\int_{\RN} u|v|^{\ss+1} dx \right| \le C_\aa \|v\|_{m(\ss+1)}^{2(\ss+1)}\:.
\end{equation*}
By (\ref{H-1}), taking $q=m(\ss+1)$, we obtain
\begin{equation*}
\|v\|_{m(\ss+1)}
\le C\|v\|^{1-\frac{m(\ss +1)-2}{2m(\ss +1)}N} \|\gd v\|^{\frac{m(\ss +1)-2}{2m(\ss +1)}N}\:,
\end{equation*}
with
\begin{equation}
\hskip-.8in
0<\frac{m(\ss+1)-2}{2m(\ss+1)}N <1\Rightarrow \ss<\frac{4}{N-2}\:.   \label{condition-3}
\end{equation}
Then (\ref{gdv}) yields
\begin{equation}
\hskip-.8in
\|\gd v(t) \|^2 \le \mathcal{H}_0 +C_\aa \|v_0\|^{2(\ss+1)-(\ss N-2)} \| \gd v\|^{\ss N-2}. \label{final-2}
\end{equation}
\\\\
For (\ref{final-1}) and (\ref{final-2}), $\|\gd v\|$ is bounded when
$\ss N-2 <2$, i.e., $\ss<\frac{4}{N}$. Therefore, the $H^1$ norm of the solution $v$ is
bounded uniformly independent of time $t$, so we can conclude that we have global solution
for any $1\le \ss <\frac{4}{N}$.
\end{proof}

In conclusion, we have shown that the Schr\"{o}dinger-Newton system (\ref{SN})
and the Schr\"{o}dinger-Helmholtz system (\ref{SNG})
admit short time unique solution when $1\le \ss< \frac{4}{N-2}$
(by definition, $\frac{4}{N-2}=\infty$ when $N\le2$) and global existence of unique solution when
$1\le \sigma < \frac{4}{N}$.  Comparing to the result of classical
NLS (\ref{NLS}) ($\ss<\frac{2}{N-2}$ for local existence and $\ss<\frac{2}{N}$ for global
existence), one expects this ``better" result
since the nonlinear terms in system (\ref{SN}) and system (\ref{SNG}) are milder than
that of the classical nonlinear Schr\"{o}dinger equation (\ref{NLS}).

\section*{Acknowledgements}
This work was supported in part by the NSF, grant no. DMS-0504619,
the BSF grant no. 2004271 and the ISF grant no. 120/06.

\end{document}